\documentclass{amsart}
\usepackage{hyperref}
\usepackage{amssymb}
\usepackage{amsrefs}

\newtheorem{theorem}{Theorem}[section]

\newtheorem{corollary}[theorem]{Corollary}
\newtheorem{lemma}[theorem]{Lemma}

\theoremstyle{remark}
\newtheorem{remark}[theorem]{Remark}

\DeclareMathOperator{\Aut}{Aut}
\DeclareMathOperator{\End}{End}
\DeclareMathOperator{\GL}{GL}
\DeclareMathOperator{\Hom}{Hom}
\DeclareMathOperator{\Jac}{Jac}
\DeclareMathOperator{\PSL}{PSL}
\DeclareMathOperator{\re}{Re}
\DeclareMathOperator{\SL}{SL}
\DeclareMathOperator{\Tr}{Tr}

\newcommand{\FF}{{\mathbf{F}}}
\newcommand{\GG}{{\mathbf{G}}}
\newcommand{\PP}{{\mathbf{P}}}
\newcommand{\QQ}{{\mathbf{Q}}}
\newcommand{\RR}{{\mathbf{R}}}
\newcommand{\ZZ}{{\mathbf{Z}}}

\newcommand\gothp{{\mathfrak{p}}}

\newcommand{\calC}{{\mathcal{C}}}

\newcommand{\emax}{E}
\newcommand{\eps}{\varepsilon}

\newcommand{\col}{\ {:}\ }        

\newcommand{\Vladut}{Vl\u adu\c t}

\newcommand\lowtilde{\hbox{\lower0.7ex\hbox{\char`\~}}}

\begin{document}

\title[Genus bounds for curves]
{Genus bounds for curves \\ with fixed Frobenius eigenvalues}

\author[Elkies]{Noam D. Elkies}    
\address{
  Department of Mathematics,
  Harvard University,
  Cambridge, MA 02138--2901
}
\email{elkies@math.harvard.edu}
\urladdr{\href{http://www.math.harvard.edu/~elkies/}
              {http://www.math.harvard.edu/\lowtilde{}elkies/}}

\author[Howe]{\hbox{Everett W. Howe}}       
\address{Center for Communications Research,
         4320 Westerra Court,
         San Diego, CA 92121-1967, USA.}
\email{however@alumni.caltech.edu}
\urladdr{\href{http://www.alumni.caltech.edu/~however/}
              {http://www.alumni.caltech.edu/\lowtilde{}however/}}

\author[Ritzenthaler]{\hbox{Christophe Ritzenthaler}}
\address{Institut de Math\'ematiques de Luminy,
         UMR 6206 du CNRS,
         Luminy, Case 907, 13288 Marseille, France.}
\email{ritzenth@iml.univ-mrs.fr}
\urladdr{\href{http://iml.univ-mrs.fr/~ritzenth/}
              {http://iml.univ-mrs.fr/\lowtilde{}ritzenth/}}
\thanks{The third author is partially supported by grant MTM2006-11391
        from the Spanish MEC and by grant ANR-09-BLAN-0020-01 from the 
        French ANR}              

\date{14 December 2011}
\keywords{Curve, Jacobian, Weil polynomial, Frobenius eigenvalue, genus,
          linear programming}

\subjclass[2010]{Primary 14G10; Secondary 11G20, 14G15, 14H25}

\begin{abstract}
For every finite collection $\calC$ of abelian varieties over~$\FF_q$,
we produce an explicit upper bound on the genus of curves over $\FF_q$
whose Jacobians are isogenous to a product of powers of elements 
of~$\calC$.

Our explicit bound is expressed in terms of the Frobenius angles of
the elements of $\calC$.  In general, suppose that $S$ is a finite 
collection of $s$ real numbers in the interval $[0,\pi]$.  If 
$S = \{0\}$ set $r = 1/2$; otherwise, let
\[
r = \#(S\cap\{\pi\}) + 2 \! \sum_{\theta\in S\setminus\{0,\pi\}}
      \left\lceil \frac{\pi}{2\theta}\right\rceil.
\]
We show that if $C$ is a curve over $\FF_q$ whose genus is greater
than
\[
\min\left(
  23\, s^2  q^{2s} \log q, 
  \
  (\sqrt{q}+1)^{2r} \left(\frac{ 1 + q^{-r}}{2} \right)
\right),
\]
then $C$ has a Frobenius angle $\theta$ such that neither $\theta$ 
nor $-\theta$ lies in $S$.

We do not claim that this genus bound is best possible.  For any 
particular set $S$ we can usually obtain a better bound by solving a
linear programming problem.  For example, we use linear programming to
give a new proof of a result of Duursma and Enjalbert:  If the Jacobian
of a curve $C$ over $\FF_2$ is isogenous to a product of elliptic curves
over $\FF_2$ then the genus of $C$ is at most~$26$.  As Duursma and
Enjalbert note, this bound is sharp, because there is an $\FF_2$-rational
model of the genus-$26$ modular curve $X(11)$ whose Jacobian splits 
completely into elliptic curves.

As an application of our results, we reprove (and correct a small error
in) a result of Yamauchi, which provides the complete list of positive
integers $N$ such that the modular Jacobian $J_0(N)$ is isogenous 
over $\QQ$ to a product of elliptic curves.
\end{abstract}

\maketitle

\section{Introduction}
\label{S:intro}

Let $(C_n)$ be a sequence of curves over a finite field $k$ such that
the genus of $C_n$ tends to infinity with $n$.  
Serre~\cite{serre}*{Cor.~2, p.~93} applies a result of Tsfasman and 
\Vladut~\cites{tsfasman,TV} to show that the dimension of the largest
$k$-simple isogeny factor of the Jacobian of $C_n$ tends to infinity 
with~$n$.  In this article we give an explicit bound for this 
asymptotic result.

The \emph{Weil polynomial} of a $d$-dimensional abelian variety $A$
over a finite field $\FF_q$ is the characteristic polynomial of the
$q$-th power Frobenius endomorphism of $A$ (acting, for instance, on
the $\ell$-adic Tate module of $A$).  The Weil polynomial is a monic
polynomial in $\ZZ[x]$ of degree $2d$, and its complex roots 
$\alpha_1,\ldots,\alpha_{2d}$ all have magnitude~$\sqrt{q}$.  The roots
can be written $\alpha_j = \sqrt{q} \exp(i \theta_j)$ for real numbers
$\theta_j$ in the half-open interval $(-\pi, \pi]$, and the roots can
be ordered so that $\alpha_j\alpha_{j+d} = q$ and
\[
0 \le \theta_1 \le \theta_2 \le \cdots \le \theta_d \le \pi.
\]
The $\theta_j$ are called the \emph{Frobenius angles} of $A$, and the
$\theta_j$ that are nonnegative are the \emph{nonnegative Frobenius
angles} of $A$.  If $C$ is a curve over $\FF_q$, we define the 
Frobenius angles of $C$ to be the Frobenius angles of its Jacobian.

We show that for any finite set $S$ of real numbers in the interval
$[0,\pi]$, every curve over $\FF_q$ of sufficiently large genus has
a nonnegative Frobenius angle that does not lie in $S$.

\begin{theorem}
\label{T:angles}
Let $S$ be a finite set of $s$ real numbers $\theta$ with 
$0\le\theta\le\pi$.  If $S = \{0\}$ set $r = 1/2$; otherwise,
take
\[
r = \#(S\cap\{\pi\}) + 2 \! \sum_{\theta\in S\setminus\{0,\pi\}}
       \left\lceil \frac{\pi}{2\theta}\right\rceil,
\]
where $\lceil x\rceil$ denotes the least integer greater than or equal
to the real number $x$.  Let 
\[
B_1 = 23\, s^2 q^{2s}\log q 
\text{\quad and \quad}
B_2 =(\sqrt{q}+1)^{2r} \left(\frac{ 1 + q^{-r}}{2}\right).
\]
If $C$ is a curve over~$\FF_q$ whose nonnegative Frobenius angles all
lie in $S$, then the genus $g$ of $C$ satisfies $g \le B_1$ and 
$g < B_2$.
\end{theorem}

\begin{remark}
Elementary calculations show that if $s > 7\sqrt{q}\log q$ 
then the bound $B_1$ from Theorem~\ref{T:angles}
is smaller than the bound $B_2$.
\end{remark}

Theorem~\ref{T:angles}
allows us to derive an explicit version of Serre's result.

\begin{corollary}
\label{C:Serre}
If $C$ is a curve of genus $g > 2$ over a finite field $\FF_q$, then
the Jacobian of $C$ has a simple isogeny factor of dimension greater
than 
\[
\sqrt{\frac{\log\log g}{6\log q}}.
\]
\end{corollary}

The first genus bound from Theorem~\ref{T:angles} leads to a corollary
that does not mention Frobenius angles.

\begin{corollary}
\label{C:varieties}
Let $A$ be a $d$-dimensional abelian variety over~$\FF_q$.  If $C$ is a
curve over $\FF_q$ of genus greater than $23\, d^2 q^{2d}\log q$ then 
the Jacobian of $C$ has a simple isogeny factor $B$ that is not an 
isogeny factor of $A$.
\end{corollary}

One natural choice for the set $S$ in Theorem~\ref{T:angles} is the
set of all nonnegative Frobenius angles for elliptic curves over a 
given finite field.  Applying the theorem to this set leads to the 
following corollary:

\begin{corollary}
\label{C:split}
Suppose $C$ is a curve over $\FF_q$ whose Jacobian is isogenous over
$\FF_q$ to a product of elliptic curves.  Then the genus of $C$ is at
most $510\, q^{8\sqrt{q} + 3} \log q$.
\end{corollary}

We prove Theorem~\ref{T:angles} and its corollaries in 
Section~\ref{S:proof}.  As we will show in Remark~\ref{R:B1}, 
for every prime power $q$ there is a family of sets $S$ for which
(in the notation of Theorem~\ref{T:angles}) the ratio $B_1/g$ is 
less than $47 s^2 \log q$; similarly, in Remark~\ref{R:B2}
we show that there is a family of examples with $s=1$ and 
$q\to\infty$ for which the ratio $B_2/g$ approaches $1$.
For many sets~$S$, however, we expect that the bounds given by 
Theorem~\ref{T:angles} are likely to be far from optimal;
in Section~\ref{S:LP} we show how tighter bounds can sometimes
be obtained by solving an integer linear programming problem. 
As a concrete example, we take $S$ to be the set of nonnegative
Frobenius angles of the elliptic  curves over $\FF_2$, and use the
linear programming method to show that the genus of a curve over
$\FF_2$ whose Jacobian splits up to isogeny into elliptic curves 
is at most~$26$, a result proved earlier by Duursma and
Enjalbert~\cite{DE} by a different (but related) method.  
In Section~\ref{S:X11} we show that for this particular set $S$
the genus bound of $26$ is sharp, because the genus-$26$ modular
curve $X(11)$ has a model over $\FF_2$ whose Jacobian is isogenous
over $\FF_2$ to a product of elliptic curves.  
(We find sharp upper bounds in other cases as well; see
Remark~\ref{R:ordinary}.)
In Section~\ref{S:X0} we use this genus bound to give a simple proof
of a result of Yamauchi~\cite{yamauchi} on the values of $N$ for which
the Jacobian of the modular curve $X_0(N)$ is isogenous over $\QQ$ to a
product of elliptic curves.  

\subsubsection*{Conventions and notation.}
Operators such as $\Hom$, $\End$ or $\Aut$ applied to varieties over a
field $k$ will always refer to $k$-rational homomorphisms and
endomorphisms. Similarly, when we say that an abelian variety over $k$
is `decomposable', `split', or `simple', these words should be
interpreted with respect to isogenies over~$k$. In the sequel, we use
$e(x)$ to denote $\exp(i x)$, where $i = \sqrt{-1}$, and we use $\re(x)$
to denote the real part of a complex number $x$.

\subsubsection*{Acknowledgments.}
After an initial version of this work was posted on the arXiv, Mike
Zieve alerted us to the paper of Duursma and Enjalbert~\cite{DE} 
mentioned above, which contains a result that implies our 
Lemma~\ref{L:key-lemma} (see~\S2.4 of the arXiv version of~\cite{DE})
and which proves that a curve over $\FF_2$ whose Jacobian is isogenous
to a product of elliptic curves has genus at most~$26$.  The print 
version of their paper also states that there is a model of $X(11)$
over $\FF_2$ whose Jacobian is isogenous over $\FF_2$ to a produce 
of elliptic curves; this is proved in an addendum added to the arXiv
version.  We are grateful to Mike Zieve for telling us about~\cite{DE},
and to Iwan Duursma for discussions about these results.

We are grateful to Armand Brumer for informing us of Yamauchi's
paper~\cite{yamauchi}.
 
In the initial version of this paper, the bound $B_1$ in 
Theorem~\ref{T:angles} was a doubly-exponential expression in~$s$ that
we obtained via an argument using our Lemma~\ref{L:key-lemma}.  We are
extremely grateful to Zeev Rudnick and Sergei Konyagin for pointing out
to us that the work of Smyth~\cite{smyth} could be used to get a bound 
that is singly-exponential in~$s$.  Our proof that $B_1$ gives an upper
bound for the genus is due almost entirely to them.

The work in this paper was begun at the GeoCrypt 2009 conference in
Pointe-\`a-Pitre.  We would like to thank the organizers of the 
conference for providing a stimulating environment for collaboration.

\section{Proof of Theorem~\ref{T:angles} and its corollaries}
\label{S:proof}

In this section we prove Theorem~\ref{T:angles} and 
Corollaries~\ref{C:varieties} and~\ref{C:split}.  The arguments that
show that $B_1$ and $B_2$ give upper bounds on the genus are 
independent of one another, so we break the proof of 
Theorem~\ref{T:angles} into two parts. 

\begin{proof}[Proof of Theorem~\textup{\ref{T:angles}}, 
              part \textup{1:} $g\le B_1$]
The bound $g\le B_1$ holds trivially when $S$ is empty, so we may 
assume that $s > 0$.  First we consider the case where $s>1$.

Suppose $C$ is a curve over $\FF_q$ all of whose Frobenius angles lie
in~$S$.  Let the elements of $S$ be $\theta_1,\ldots,\theta_s$, and 
for each index $j$ let $b_j$ denote the multiplicity of $\theta_j$ as
a Frobenius angle for $C$ (unless $\theta_j = 0$ or $\theta_j = \pi$,
in which case we let $b_j$ denote half of this multiplicity), so that
we have
\[ 
b_1 + \ldots + b_s = g. 
\]
Reindex the $\theta_j$ and the corresponding $b_j$ so that 
$b_1\ge b_2 \ge \cdots \ge b_s$, and note that this implies that
$b_1 \ge g / s$. Finally, for every $j$ let $z_j = e(\theta_j)$, and
for every integer $k > 0$ let $G(k)$ denote the weighted power sum
\[
G(k) = b_1 z_1^k + \cdots + b_s z_s^k. 
\]
Weil's `Riemann Hypothesis' for curves over finite fields says that
\[
\#C(\FF_{q^k}) = q^k + 1 - \sum_{j=1}^{2g} \alpha_j^k
              = q^k + 1 - 2 q^{k/2} \re G(k),
\]
so we find that
\[ 
2 \re G(k) = q^{k/2} + q^{-k/2} - \#C(\FF_{q^k}) / q^{k/2}.
\]

Now we apply a result of Smyth~\cite{smyth}. Let
\[ 
\beta := (b_2 + \cdots + b_s)/b_1,
\]
let $\lambda$ be an arbitrary real number with $0 < \lambda \le 1$,
and let 
\[
K := 1 + \lfloor(4\beta+3)/\lambda\rfloor.
\]
Smyth shows that then
\[ 
(b_1/4)(1-\lambda) < \max_{1\le k \le K} \re G(k).
\]
In our case, we have $\beta = g/b_1 - 1 \le s - 1$, so 
$K \le 1 + (4\beta+3)/\lambda  \le L$, where we set
\[ 
L:= 1 + (4s - 1)/\lambda.
\]
Also, for each $k$ we have $\re G(k) \le (1/2)(q^{k/2} + q^{-k/2})$.
Since $q^x + q^{-x}$ is an increasing function of $x$ for $x>0$, 
Smyth's result shows that 
\[ 
(b_1/4)(1-\lambda) < (1/2)(q^{L/2} + q^{-L/2}).
\]
As we noted earlier, we have $b_1\ge g/s$, so for every $\lambda$ 
with $0<\lambda<1$ we have
\begin{equation}
\label{EQ:Smyth}
\frac{g}{s} < \frac{2}{1-\lambda}\left(q^{L/2} + q^{-L/2}\right).
\end{equation}

We will apply this inequality with $\lambda = 1 - 1/a$, where
\[ 
a = \left(\frac{4s-1}{2}\right)\log q. 
\]
(This choice of $\lambda$ was obtained by taking the first two terms of
the power series expansion (in $1/a$) of the value of $\lambda$ that
minimizes $q^{L/2}/(1-\lambda)$.)  Note that $L > 4s \ge 8$, so that
\[
q^{L/2} + q^{-L/2} = q^{L/2}(1 + q^{-L}) < \frac{257}{256}q^{L/2}.
\]
Combining this with~\eqref{EQ:Smyth}, we find
\begin{align*}
\frac{g}{s} &< \left(\frac{257}{256}\right)
               \left(\frac{2}{1-\lambda}\right) q^{L/2} \\
            &= \left(\frac{257}{256}\right)
               (4s-1)(\log q) \, \exp((L/2)\log q).
\end{align*}
Using the equality $L = 1 + (4s - 1)/\lambda$, we find that
\begin{align*}
(L/2) \log q    &= (\log q)/2 + \left(\frac{4s - 1}{2}\right)
                                      \frac{\log q}{\lambda} \\
                &= (\log q)/2 + a/\lambda \\
                &= (\log q)/2 + a^2/(a-1) \\
                &= (\log q)/2 + a + a/(a-1)\\
                &= 2 s \log q + a/(a-1).
\end{align*}
Since $s\ge 2$ and $q\ge 2$ we have $a\ge (7/2)\log 2$ and
$\exp(a/(a-1)) < 5.5$.  Thus
\[ 
\exp((L/2)\log q) < 5.5\, \exp(2s\log q) = 5.5 \,q^{2s},
\]
so
\[
\frac{g}{s} < \left(\frac{257}{256}\right)
              (4s)(\log q) (5.5\, q^{2s})
            < 23\, s q^{2s} \log q.
\]
This shows that $g < B_1$ when $s > 1$.

Now suppose $s = 1$.  Let $\theta$ be the unique element of $S$.
If $\theta \ge \pi/2$ then the quantity $r$ defined in the statement 
of the theorem is equal to $2$.  It is easy to show that then 
$B_2 < B_1$, so the bound $g < B_2$, proved below, shows that 
$g < B_1$.

The final case to consider is when $s = 1$ and $\theta < \pi/2$.  
If $C$ has only one nonnegative Frobenius angle $\theta$, then 
$e(\theta)\sqrt{q}$ must be an algebraic integer of degree at most
$2$ over the rationals, so the quantity 
$t =  e(\theta)\sqrt{q} +  e(-\theta)\sqrt{q}$ is an integer.
Furthermore, since $\theta$ is less than $\pi/2$ the integer $t$ 
is positive.  Then Weil's theorem shows that 
\[
0 \le \#C(\FF_q) = q + 1 - gt \le q + 1 - g,
\]
so that $g \le q + 1.$  Again, it follows easily that $g < B_1$.
\end{proof}

Our proof that $B_2$ gives an upper bound for the genus relies upon
the following lemma.

\begin{lemma}
\label{L:key-lemma}
Let $S$ be a set of real numbers in $[0,\pi]$ and let $T \in \RR[X]$ be
a polynomial with nonnegative coefficients such that $T(0)=0$ and
$\re(T(e(\theta))) \geq 1$ for all $\theta \in S$.  If $C$ is a curve 
over $\FF_q$ whose nonnegative Frobenius angles all lie in~$S$, then 
the genus of $C$ is at most $(T(q^{1/2}) + T(q^{-1/2})) / 2$.
\end{lemma}

\begin{proof}
Let $\alpha_1,\ldots,\alpha_{2g}$ be the complex roots of the Weil 
polynomial of~$C$, listed with appropriate multiplicities, so that
from Weil's theorem we have
\[
\#C(\FF_{q^m}) = q^m + 1 - \sum_{j=1}^{2g} \alpha_j^m
\]
for all integers $m>0$.  Write 
\[
T = a_1 x + \cdots + a_n x^n, \quad a_j \geq 0.
\]
Then 
\begin{align*}
0 &\le \sum_{m=1}^n a_m \frac{\#C(\FF_{q^m}) }{q^{m/2}}\\
   & =  \sum_{m=1}^n a_m \left(
         q^{m/2} + q^{-m/2} - \sum_{j=1}^{2g} \alpha_j^m/q^{m/2}
         \right)\\
   &= T(q^{1/2}) + T(q^{-1/2}) - \sum_{j=1}^{2g} T(\alpha_j/\sqrt{q}).
\end{align*}
By hypothesis, each summand $T(\alpha_j/\sqrt{q})$ has real part at 
least $1$, so 
\[ 
0 \le T(q^{1/2}) + T(q^{-1/2}) - 2g,
\]
and hence 
\[ 
g \le \frac{T(q^{1/2}) + T(q^{-1/2}) }{2}.
\]
\end{proof}

\begin{remark}
This lemma could also be proved by using the results in~\S2.4 of the
arXiv version of~\cite{DE}.
\end{remark}

\begin{remark}
Lemma~\ref{L:key-lemma} is stated for polynomials $T$ with nonnegative
coefficients, but an analogous statement holds when $T$ is a power 
series in $\RR[[x]]$ with nonnegative coefficients, so long as its 
radius of convergence exceeds $\sqrt{q}$.
\end{remark}

With this lemma in hand, we complete the proof of 
Theorem~\ref{T:angles}.  Our proof will depend on a careful choice of
the polynomial $T$.

\begin{proof}[Proof of Theorem~\textup{\ref{T:angles}}, 
              part \textup{2:} $g < B_2$]
It is easy to check that $g < B_2$ when $S$ is empty, so we may assume
that $S$ is nonempty. If $S = \{0\}$ then the quantity $r$ from the
theorem is equal to $1/2$. In this case we can take $T = x$, and we find
that $(T(q^{1/2}) + T(q^{-1/2}))/2$ is less than $B_2$. Then 
Lemma~\ref{L:key-lemma} shows that $g < B_2$.  So 
let us assume that $S$ contains a nonzero element.

Given a nonzero $\theta\in S$, let
$m = \left\lceil \frac{\pi}{2\theta}\right\rceil,$  so that
$\cos(m\theta)\le 0$.  If $\theta = \pi$ let $P_\theta$ be the
polynomial $1+x$; otherwise, set 
$P_\theta = 1 - 2\cos(m\theta)x^m + x^{2m}$. In both cases we have
$P_\theta(e(\theta)) = 0.$  Let 
\[
P = \prod_{\theta\in S\setminus\{0\}} P_\theta.
\]
Then $P$ is a polynomial with constant term $1$, with nonnegative 
coefficients, and with degree equal to $r$, where
\[
r = \#(S\cap\{\pi\}) + 2\sum_{\theta\in S\setminus\{0,\pi\}}
      \left\lceil \frac{\pi}{2\theta}\right\rceil
\]
is as in the statement of the theorem.  Let $T = (P-1)^2$.  Then 
$T(0) = 0$, the coefficients of $T$ are nonnegative, and for every
$\theta \in S\setminus\{0\}$ we have $T(e(\theta)) = 1.$  Also, for
each $\theta$ we have $P_{\theta}(1)\ge 2$, so $P(1)\ge 2$ and 
$T(1)\ge 1$.  Lemma~\ref{L:key-lemma} shows that any curve whose 
nonnegative Frobenius angles all lie in $S$ must have genus no larger 
than $(T(q^{1/2}) + T(q^{-1/2}))/2$.  Now, for positive real numbers 
$z$ we have $1\le P(z) \le (1 + z)^r$, so that $T(z) < (1+z)^{2r}.$  
Thus we have
\begin{align*}
T(q^{1/2}) + T(q^{-1/2}) &< (\sqrt{q} + 1)^{2r} + (1/\sqrt{q} + 1)^{2r}\\
                         &= (\sqrt{q} + 1)^{2r }(1 + q^{-r}),
\end{align*}
which gives the inequality $g < B_2$.
\end{proof}

\begin{remark}
\label{R:B2}
In Section~\ref{S:LP} we will see that our bounds can sometimes be bad.
However, the following easy example shows that at least in one case our
second bound is asymptotically exact as $q \rightarrow \infty$.  For 
any prime power~$q$ let $E$ be a supersingular elliptic over $\FF_q$
with Weil polynomial $x^2+q$, corresponding to the set $S=\{\pi/2\}$.
Applying Theorem~\ref{T:angles}, we see that the genus of a curve
$C / \FF_q$ with Jacobian isogenous to a power of~$E$ is bounded above 
by
\[
(\sqrt{q}+1)^4 \cdot \left(\frac{1+q^{-2}}{2}\right) \sim q^2/2
                  \text{\qquad as $q\to\infty$}.
\]
On the other hand, $\Jac(C)$ is isogenous to a power of~$E$ if and only
if $C$ is optimal over $\FF_{q^2}$, in the sense that its number of 
points attains the Weil upper bound.  But the maximal genus of such a 
curve is $q(q-1)/2 \sim q^2/2$ (see \cite{ihara}), attained by the
Hermitian curve $H_q$ defined by $x^{q+1}+y^{q+1}+z^{q+1}=0$ (see for 
instance~\cite{RS}).
\end{remark}

\begin{remark}
\label{R:moreB2}
The Hermitian curve $H_q$, viewed as a curve over $\FF_{q^2}$, again
gives an example of a curve whose genus comes close to the upper
bound $B_2$.  If we take $S = \{\pi\}$ then the bound $B_2$ for the
field $\FF_{q^2}$ is 
\[
(q+1)^2 \cdot \left(\frac{1+q^{-2}}{2}\right) \sim q^2/2
                  \text{\qquad as $q\to\infty$},
\]
while $H_q$ has genus $q(q-1)/2$.
\end{remark}

\begin{remark}
\label{R:B1}
Hermitian curves can also be used to give examples that limit the 
extent to which we might hope to improve the bound~$B_1$.  For any 
integer $s>0$, let $S$ be the $s$-element set
\[
S = \left\{\frac{\pi}{2s}, \frac{3\pi}{2s}, \ldots, \frac{(2s-1)\pi}{2s} \right\}.
\]
Let $q$ be any prime power and set $Q=q^s$.  
Note that the nonnegative Frobenius angles
of a curve $C$ over $\FF_q$ are contained in $S$ if and only if the 
only nonnegative Frobenius angle of the base extension of $C$ to 
$\FF_Q$ is $\pi/2$.  Let $C$ be the curve
$x^{Q+1}+y^{Q+1}+z^{Q+1}=0$ over $\FF_q$.  As noted in 
Remark~\ref{R:B2}, the base extension of $C$ to $\FF_Q$ has $\pi/2$
as its unique Frobenius angle, so all of the Frobenius angles of $C$
itself are contained in $S$.  Since the genus of $C$ is 
$(q^{2s}-q^s)/2$, we see that Theorem~\ref{T:angles} would be false if
we replaced $B_1$ with any expression of the form
\[
(\text{polynomial in $s$ and $\log q$}) \, q^{2s-\eps}
\]
for a positive constant~$\eps$.
\end{remark}

We end this section by proving the corollaries from the introduction.

\begin{proof}[Proof of Corollary~\textup{\ref{C:Serre}}]
We begin by noting that for every integer $n>0$, 
the number of isogeny classes of $n$-dimensional
abelian varieties over~$\FF_q$ is less than $6^{n^2}q^{n(n+1)/4}$.  
This is easy to check for $n=1$ (see the proof of 
Corollary~\ref{C:split} below) and for $n=2$, while for $n > 2$ it
follows from the results in~\cite{DH1} (as corrected in~\cite{DH2}).
We leave the details of the argument to the reader, but we do at least
note that it is helpful to observe that for every $n$, the quantity 
$v_n$ in~\cite{DH1}*{Thm.~1.2} is bounded above by $264$.
It follows that the number of isogeny classes of \emph{simple} abelian 
varieties over $\FF_q$ of dimension \emph{at most} $n$ is also less than
$6^{n^2}q^{n(n+1)/4}$, and that the number of nonnegative Frobenius
angles of abelian varieties over $\FF_q$ of dimension at most $n$ is
less than $n 6^{n^2}q^{n(n+1)/4}$.

Suppose, to obtain a contradiction, that the corollary is false, and
let $C$ be a curve of genus $g$ over $\FF_q$ that provides a 
counterexample.  Take
\[
n = \left\lfloor\sqrt{\frac{\log \log g}{6 \log q}}\right\rfloor.
\]
For $C$ to provide a counterexample we must have~$n\ge 1$.
We will apply Theorem~\ref{T:angles} to the set of nonnegative 
Frobenius angles of abelian varieties over $\FF_q$ of dimension at
most~$n$; as we have just noted, this means that 
in Theorem~\ref{T:angles} we have $s < n 6^{n^2} q^{n(n+1)/4}$. 
Simple estimates for the terms
in the right-hand side of this inequality show that 
\begin{equation}
\label{EQ:s}
s < q^{4n^2}.
\end{equation}
Note that there are at least $5$ isogeny classes of elliptic curves
over any field, so $s\ge 5$; using this fact it is easy to show 
that the bound $B_1 = 23 s^2 q^{2s}\log q$ from the theorem satisfies
\begin{equation}
\label{EQ:B1}
B_1 < q^{4s}.
\end{equation}

The definition of $n$ tells us that 
\[
\log\log g \ge 6 n^2 \log q,
\]
so exponentiating gives us
\[
\log g \ge q^{6n^2} > 4 q^{4n^2} \log q > 4 s \log q,
\]
where the third inequality follows from~\eqref{EQ:s}.  Exponentiating
again, we find that
\[ 
g > q^{4s} > B_1,
\]
by~\eqref{EQ:B1}.  Theorem~\ref{T:angles} then shows that $\Jac C$ has
a simple isogeny factor of dimension greater than $n$, contradicting 
our assumption that $C$ was a counterexample to the corollary.  This 
contradiction completes the proof.
\end{proof}

\begin{proof}[Proof of Corollary~\textup{\ref{C:varieties}}]
Let $S$ be the set of nonnegative Frobenius angles of the abelian 
variety $A$ over $\FF_q$, so that $\#S\le d$.  If $C$ is a curve over
$\FF_q$ with genus greater than $23\, d^2 q^{2d}\log q$ then $C$ has 
a nonnegative Frobenius angle $\theta$ not in $S$.  Let $B$ be any 
element of the unique isogeny class of simple abelian varieties over
$\FF_q$ that have $\theta$ as one of their Frobenius angles.  Then 
$B$ is an isogeny factor of $\Jac C$.  On the other hand, $A$ and 
$B$ have coprime Weil polynomials, so by the Honda--Tate theorem $B$
is not an isogeny factor of~$A$.
\end{proof}

\begin{proof}[Proof of Corollary~\textup{\ref{C:split}}]
Let $m$ be the largest integer with $m^2 \le 4q$.  The Weil polynomial
of an elliptic curve over $\FF_q$ is of the form $x^2 - tx + q$, where
$t$ is an integer with $-m\le t\le m$.  (Not every $t$ in this range 
need come from an elliptic curve.)  The nonnegative Frobenius angle 
$\theta_t$ corresponding to a given $t$ is 
$\theta_t = \cos^{-1} (t / (2 \sqrt{q} )).$  We take $S$ to be the
set $\{\theta_t \col -m \le t\le m\}$.

Clearly $\#S \le 4\sqrt{q} + 1$.  The bound $B_1$ from 
Theorem~\ref{T:angles} is then
\begin{align*}
23\, s^2 q^{2s} \log q &\le 23\, (4\sqrt{q} + 1)^2 q^{8\sqrt{q} + 2} \log q\\
                       &\le 23\, \left(\frac{4\sqrt{2} + 1}{\sqrt{2}}\right)^2
                          q^{8\sqrt{q} + 3} \log q \\
                       &< 510\, q^{8\sqrt{q} + 3} \log q.
\end{align*}                     
\end{proof}

\section{Linear programming}
\label{S:LP}

The arguments we used to prove Theorem~\ref{T:angles} used no special
properties of the set $S$ and give a bound that is almost certainly
far from optimal when applied to most~$S$.  For example, suppose we
take $S$ to be the set of nonnegative Frobenius angles coming from 
the elliptic curves over~$\FF_2$.  There are five isogeny classes of
elliptic curves over~$\FF_2$ (each containing a single curve),
corresponding to the five possible traces
of Frobenius $-2,-1,0,1,2$ (see \cite{deuring}, 
\cite{waterhouse}*{Theorem 4.1}). For each of these traces~$t$, let
$\alpha_t = (-t + \sqrt{t^2-8})/2$, let $\theta_t$ be the argument of 
$\alpha_t$, and let $E_t$ be an elliptic curve over $\FF_2$ with
trace~$t$.  Let us apply Theorem~\ref{T:angles} to the set 
$S =\{\theta_t\}$.  We find that $B_1 \approx 408125$ and 
$B_2 \approx 2.6\times10^{10}$.  However, the best upper bound on the 
genus is much smaller than these numbers.

\begin{theorem}[Duursma and Enjalbert~\cite{DE}]
\label{T:F2}
Suppose $C$ is a curve over $\FF_2$ whose Jacobian is isogenous to a 
product of powers of the $E_t$.  Then the genus of $C$ is at most~$26$.
\end{theorem}

\begin{remark}
Duursma and Enjalbert prove this result by using a stronger version of
our Lemma~\ref{L:key-lemma}.  We give the proof presented here as an 
example of our general technique of using linear programming to get
genus bounds.
\end{remark}

\begin{remark}
We will see in Section~\ref{S:X11} that the bound in Theorem~\ref{T:F2}
is sharp.
\end{remark}

\begin{proof}[Proof of Theorem~\textup{\ref{T:F2}}]
Our bounds in Theorem~\ref{T:angles} were obtained from the fact that
the number of points on a curve over a finite field is always 
nonnegative.  We actually know a somewhat stronger constraint:  For 
every $n>0$, the number of degree-$n$ places on a curve is
nonnegative.  If the size of the base field is large compared to the 
genus of the curve in question, the bounds we get from using place 
counts are not much better than the one we get from point counts.  
However, Theorem~\ref{T:F2} involves a very small field indeed.

Let $C$ be a curve as in the statement of the theorem.  In particular,
suppose there are integers $e_t \ge 0$ such that the Jacobian of $C$
is isogenous to
\[
E_{-2}^{e_{-2}} \times
E_{-1}^{e_{-1}} \times
E_{ 0}^{e_{ 0}} \times
E_{ 1}^{e_{ 1}} \times
E_{ 2}^{e_{ 2}}.
\]
Then for every $n>0$ the number of points on $C$ over $\FF_{2^n}$ is
given by
\[
\#C(\FF_{2^n}) = 2^n + 1 - \sum_{-2\le t\le 2} e_t\Tr(\alpha_t^n),
\]
and the number $N_n$ of degree-$n$ places on $C$ is given by
\[
N_n = \frac{1}{n}\sum_{d|n} \mu(n/d)\#C(\FF_{2^d}),
\]
where $\mu$ is the M\"obius function.

The nonnegativity of the number of degree-$n$ places, for 
$n=1,\ldots,8$, is expressed by the following inequalities:
\begin{equation}
\label{EQ:places8}
\begin{array}{rrrrrl}
- 2 e_{-2} & -   e_{-1} &          & +   e_1 & + 2 e_2 & \le  3  \\
    e_{-2} & -   e_{-1} & - 2 e_0  & - 2 e_1 & -   e_2 & \le  1  \\
  2 e_{-2} & + 2 e_{-1} &          & - 2 e_1 & - 2 e_2 & \le  2  \\
- 2 e_{-2} & +   e_{-1} & + 3 e_0  & +   e_1 & - 2 e_2 & \le  3  \\
  2 e_{-2} & - 2 e_{-1} &          & + 2 e_1 & - 2 e_2 & \le  6  \\
-   e_{-2} & +   e_{-1} & - 2 e_0  & + 3 e_1 & +   e_2 & \le  9  \\
- 2 e_{-2} & + 2 e_{-1} &          & - 2 e_1 & + 2 e_2 & \le 18  \\
  5 e_{-2} & - 4 e_{-1} & + 3 e_0  & - 4 e_1 & + 5 e_2 & \le 30  \\
\end{array}
\end{equation}
We claim that $\sum_{t=-2}^2 e_t \leq 26$ for any nonnegative integers
$e_t$ satisfying the inequalities above.  Indeed we show this even if
the $e_t$ are allowed to be nonnegative real numbers.  Maximizing their
sum is then a linear programming problem.  Solving the dual problem, we
find that if we take $39$ times the third inequality 
of~\eqref{EQ:places8}, plus $44$ times the fourth inequality, plus $78$
times the sixth inequality, plus $32$ times the eighth inequality, we 
obtain
\[
 72 e_{-2}  + 72 e_{-1}  + 72 e_0   + 72 e_1  + 72 e_2  \le 1872,
\]
so that $e_{-2}+e_{-1}+e_0+e_1+e_2\le 1872/72 = 26$, as claimed.
\end{proof}

\begin{remark}
\label{R:genus26}
By enumerating $5$-tuples of nonnegative integers
$(e_{-2},e_{-1},e_0,e_1,e_2)$ that sum to $26$ and checking whether
they satisfy the inequalities~\eqref{EQ:places8}, one finds that if a
genus-$26$ curve over $\FF_2$ has a completely split Jacobian, then 
the Jacobian is isogenous to one of the varieties
\begin{align*}
&E_{-2}^4 \times E_{-1}^7 \times E_0^5 \times E_1^3 \times E_2^7,\\
&E_{-2}^5 \times E_{-1}^6 \times E_0^5 \times E_1^4 \times E_2^6,
\text{\ or}\\
&E_{-2}^6 \times E_{-1}^5 \times E_0^5 \times E_1^5 \times E_2^5.
\end{align*}
\end{remark}

\begin{remark} 
We would not get a genus bound if we looked only at the inequalities
coming from counting the places of degree strictly less than $8$.  The
reader can check that for any even $g$, the values $e_{-2} = e_2 = g/2$
and $e_{-1}=e_0=e_1=0$ satisfy the first seven inequalities in the 
system~\eqref{EQ:places8}. On the other hand, \emph{a priori} there is
no reason to think that the bound we get from looking at the counts of
places of degree $8$ and less will be the best possible; perhaps by 
using the nonnegativity of several more place counts, we would get a 
better bound.  In this particular instance, results of
Section~\ref{S:X11} show that the bound we obtain is in fact sharp.
\end{remark}

\begin{remark}
A similar argument, using places of degree at most $12$, shows that a 
curve over $\FF_3$ whose Jacobian is isogenous to a product of elliptic
curves must have genus less than or equal to $2091$.   By extending the
argument slightly (or by using an integer linear programming package, 
such as the one in Magma~\cite{magma}) one can improve this upper bound
to $2085$.  We suspect that this bound is not sharp!
\end{remark}

\begin{remark}
\label{R:ordinary}
If one restricts to curves whose Jacobians are isogenous to a product 
of powers of ordinary elliptic curves, one finds a genus bound of
$3$ over $\FF_2$ and $26$ over $\FF_3$. The first bound is
reached by the the curve
\[ 
x^4 + y^4 + z^4 + x^2 y^2 + x^2 z^2 + y^2 z^2 + x^2 y z + x y^2 z + x y z^2 = 0
\]
which is a twist of the Klein quartic;  the second bound is reached as
well, as we show in Section~\ref{S:X11}.
\end{remark}

\section{The modular curve $X(11)$}
\label{S:X11}

In this section we show that the modular curve $X(11)$ has a model
defined over~$\FF_2$ whose Jacobian is isogenous over $\FF_2$ to
a product of elliptic curves. Since $X(11)$ has
genus~$26$, this example shows that the bound of Theorem~\ref{T:F2} is 
sharp.  Duursma and Enjalbert~\cite{DE} provide a different proof 
that $X(11)$ has a model over $\FF_2$ with completely split Jacobian; 
their argument, found in the addendum to the arXiv version of their 
paper, relies on working with Klein's explicit model of $X(11)$ in~$\PP^4$.

Let $G$ be the twist of the finite group scheme $\ZZ/11\ZZ$ over 
$\FF_2$ on which the $\FF_2$-Frobenius acts as multiplication by~$3$, 
and let $G'$ be the Cartier dual of $G$, so that $G'$ is the twist of 
$\ZZ/11\ZZ$ on which the $\FF_2$-Frobenius acts as multiplication 
by~$2/3=-3 \bmod 11$.  Let $e\col G\times G'\to\GG_m$ be the natural
pairing from $G\times G'$ to the multiplicative group.

Let $X$ be the modular curve over $\FF_2$ whose non-cuspidal
$K$-rational points, for every extension $K$ of $\FF_2$, parametrize 
pairs $(E, \varphi)$, where $E$ is an elliptic curve over $K$ and 
$\varphi$ is an isomorphism from the group scheme $E[11]$ to the group
scheme $(G\times G')\otimes_{\FF_2} K$ that takes the Weil pairing on
$E[11]$ to the pairing~$e$.

\begin{theorem}
\label{T:modular}
The genus of $X$ is $26$, and the Jacobian of $X$ is isogenous to
\[
E_{-2}^6 \times E_{-1}^5 \times E_0^5 \times E_1^5 \times E_2^5,
\]
where each $E_t$ is an elliptic curve over $\FF_2$ with trace $t$.
\end{theorem}

\begin{proof}
The curve $X$ is geometrically isomorphic to $X(11)$, so it has
genus~$26$ and geometric automorphism group isomorphic to~$\PSL_2(11)$
(see \cite{ritzenthaler}*{Th\'eor\`eme~5} and~\cite{BCG}).  Consider
the group scheme $(G\times G')\otimes_{\FF_2}\FF_4$; it is simply 
$(\ZZ/11\ZZ)^2$, with the $\FF_4$-Frobenius acting as multiplication
by~$-2$.  The automorphism group of this $\FF_4$-scheme is $\GL_2(11)$,
and the subgroup of automorphisms that respect the pairing $e$ is
isomorphic to $\SL_2(11)$.  There is a surjective map
\[
\Aut \left( (G\times G')\otimes_{\FF_2}\FF_4, e\right)
\to \Aut (X\otimes_{\FF_2}{\FF_4})
\]
that sends an automorphism $\alpha$ of the finite group-scheme to the
automorphism $\beta$ of $X\otimes_{\FF_2}\FF_4$ that takes a pair 
$(E,\varphi)$ to $(E,\alpha\varphi)$, and the kernel of this map is the
group $\{\pm1\}$.  Therefore all of the geometric automorphisms of $X$ 
are already defined over~$\FF_4$.

Using~\cite{GGHZ}*{Lemma~2.1} we see that the twists of the curve
$X\otimes_{\FF_2}\FF_4$ correspond to the conjugacy classes of 
$\Aut (X\otimes_{\FF_2}\FF_4)$; also, the automorphism group of the
twist corresponding to the conjugacy class of an element $\alpha$ is
isomorphic to the commutator of $\alpha$.  Since $\PSL_2(11)$ has 
trivial center, we see that every nontrivial twist of 
$X\otimes_{\FF_2}\FF_4$ has automorphism group strictly smaller 
than $\PSL_2(11)$.

Now take the $\QQ(\sqrt{-11})$-rational model $Y$ of $X(11)$ considered
by Ligozat~\cite{ligozat}*{Example~3.7.3, pp.~199--200}; Ligozat calls
this curve $X(11)_{K_{11}}$.  Let $\gothp$ be the prime of 
$\QQ(\sqrt{-11})$ over~$2$, with residue field $\FF_4$.  The 
automorphism group of $Y$ is $\PSL_2(11)$, so the reduction of $Y$ 
modulo $\gothp$ also has automorphism group $\PSL_2(11)$.  Therefore,
the reduction of $Y$ must be $X\otimes_{\FF_2}\FF_4$.  Applying a 
result of Ligozat~\cite{ligozat}*{Prop 3.6.1, p.~223}, we find that the
Jacobian of $X$ is isogenous to 
$E_{-2}^{11}\times E_0^5 \times E_1^{10}$.  (To see this, we must take
the elliptic curves mentioned in Ligozat's proposition and compute 
their reductions modulo $\gothp$.)

It follows that over $\FF_2$, the Jacobian of $X$ is isogenous to
\[
E_{-2}^a E_{-1}^b E_0^5 E_1^{10-b} E_2^{11-a}
\]
for some choice of $a$ and $b$, and by Remark~\ref{R:genus26}, we know
that $(a,b)$ is one of $(4,7)$, $(5,6)$, and~$(6,5)$.

For each of the three possible isogeny classes we can compute the
associated zeta function and the number of $\FF_2$-rational points on
a curve with that zeta function.   We find that if the number of
$\FF_2$-rational points on $X$ is $1$, then $(a,b)=(4,7)$; if
$\#X(\FF_2)=3$ then $(a,b) = (5,6)$; and if $\#X(\FF_2) = 5$ then
$(a,b) = (6,5)$.

Consider $E_0$, the unique elliptic curve over $\FF_2$ with trace~$0$.
The characteristic polynomial of Frobenius on $E_0[11]$ is $x^2 + 2$,
so the action of Frobenius on $E_0[11]$ has two eigenspaces, one with
eigenvalue $3$ and one with eigenvalue $-3$.  One finds that there are
$10$ pairing-respecting isomorphisms $\varphi:E_0[11]\to G\times G'$, 
and since $(E_0,\varphi)$ and $(E_0,-\varphi)$ are represented by the 
same point on $X$, we have found $5$ $\FF_2$-rational points on $X$.
Therefore the Jacobian of $X$ decomposes as claimed in the statement 
of the theorem.
\end{proof}

\begin{remark}
One can show that the curve $X$ has $60$ cusps (that is, points that 
lie over the point at infinity on the the $j$-line).  The field of 
definition of $10$ of the cusps is $\FF_{2^5}$; the field of definition
of the other $50$ is $\FF_{2^{10}}$.  Using these facts, together with 
the modular interpretation of the non-cuspidal points on $X$, we can 
compute the number of points on $X$ over any (reasonably small) 
extension of~$\FF_2$. This gives another method of computing the 
decomposition of the Jacobian of~$X$.
\end{remark}

\begin{remark}
Applying Ligozat's Proposition~3.6.1, we see that the Jacobian of
$X(11)_{K_{11}} \otimes \FF_3$ splits into a product of ordinary 
elliptic curves. Hence the bound $26$ for ordinary elliptic curves over
$\FF_3$ from Remark~\ref{R:ordinary} is reached as well.
\end{remark}

\section{Application: Modular curves with split Jacobians}
\label{S:X0}

Let $J_0(N)$ denote the Jacobian of the modular curve $X_0(N)$
over~$\QQ$.  Cohen~\cite{cohen} (mentioned
in~\cite{serre}*{Remarque~2, p.~90} with the value $N=27$ omitted)
has computed a list of the odd integers $N$ for which $J_0(N)$ is 
isogenous to a product of elliptic curves, and 
Yamauchi~\cite{yamauchi}*{Thm.~1.1}
extended this list to include even values of $N$ as well.  In this section we use 
Theorem~\ref{T:F2} and the mathematical software package 
Sage~\cite{sage} to recompute Yamauchi's list; we note that the list
in Yamauchi's theorem mistakenly includes $N=672$ and omits $N=28$.

\begin{theorem}
\label{T:X0}
Let $N = 2^e n$ with $n$ odd.  Then $J_0(N)$ is isogenous over~$\QQ$
to a product of elliptic curves if and only if $n$ appears in the
following table and $e \leq \emax(n)$ with $\emax(n)$ as 
tabulated\textup{:}
\[
\begin{array}{r|ccr|ccr|c}
 n & \emax(n) &\qquad &  n & \emax(n) &\qquad &   n & \emax(n) \\
\cline{1-2}\cline{4-5}\cline{7-8}
 1 & 7    & & 15 & 4    & &  37 & 0    \\
 3 & 7    & & 17 & 1    & &  45 & 4    \\
 5 & 4    & & 19 & 2    & &  49 & 0    \\
 7 & 4    & & 21 & 4    & &  57 & 1    \\
 9 & 7    & & 25 & 4    & &  75 & 4    \\
11 & 2    & & 27 & 4    & &  99 & 2    \\
13 & 2    & & 33 & 2    & & 121 & 0    \\
\end{array}
\]
\end{theorem}



\begin{proof}
Suppose $N$ is an odd integer such that $J_0(N)$ splits into elliptic
curves.  Since $N$ is odd the modular curve $X_0(N)$ has good reduction
modulo~$2$, and the reduced curve over $\FF_2$ has split Jacobian.  By 
Theorem~\ref{T:F2}, the genus of $X_0(N)$ is at most $26$.  Using the 
fact~\cite{CWZ} that the genus of $X_0(N)$ is greater than or equal to
$(N-5\sqrt{N}-8)/12$, we find that $N$ is at most~$422$.

Using the command \verb|J0(N).decomposition()| of the mathematics 
package Sage for all \verb|N| less than $423$, we find that the odd
values of~$N$ with $J_0(N)$ split are the twenty-one odd integers that 
appear as $n$ in the table.

To complete the proof, we note that if $J_0(N)$ is split then so is
$J_0(n)$ for every divisor $n$ of $N$.  Therefore the integers we are
searching for can be written $2^e n$ for some exponent $e$ and for some
$n$ among the odd values that we have just computed.  For each of the 
possible odd parts $n$, we use Sage to compute the decomposition of 
$J_0(2^e n)$ for increasing values of $e$ until we reach a Jacobian 
that does not split.  (In practice, we do not have to compute the 
decomposition of $J_0(2^e n)$ if we already know that $J_0(2^e m)$
does not split for some divisor $m$ of~$n$.  For example, since 
$J_0(2^8)$ does not split we must have $\emax(n)<8$ for each $n$.)
The largest value of~$e$ for which $J_0(2^e n)$ splits is recorded
in the table as $\emax(n)$.
\end{proof}

\begin{remark}
Ekedahl and Serre~\cite{ES} give a list of various values of $g$ such
that there exists a curve of genus $g$ over $\QQ$ with a completely 
split Jacobian.  Many of their values of $g$ come from modular 
curves~$X_0(N)$.  The tables of modular forms that they had access to
did not include values of $N$ greater than $1000$, so they missed a few
of the values from Theorem~\ref{T:X0}.  The modular curves $X_0(1152)$ 
and $X_0(1200)$, of genus $161$ and $205$, respectively, allow us to 
add two more values of $g$ to their 
list~\cite{ES}*{Th\'eor\`eme, p.~509}.

Ekedahl and Serre also seem to have missed the curve $X_0(396)$ of
genus $61$, but they obtain a curve of genus $61$ with split Jacobian
by considering a quotient of $X_0(720)$ by an involution.
\end{remark}







\end{document}